\documentclass[12pt]{article}
\usepackage{amsmath}
\usepackage{amssymb}
\usepackage{vatola}

\def\q{\quad}

\def\mod{\pmod}
\def\t{\text}
\def\f{\frac}
\def\e{\equiv}

\def\phq#1{\varphi(q^{#1})}
\def\psq#1{\psi(q^{#1})}
\def\qtq#1{\q\t{#1}\q}
\def\sls#1#2{(\f{#1}{#2})}
 
\def\Ls#1#2{\Big(\f{#1}{#2}\Big)}
\let \pro=\proclaim
\let \endpro=\endproclaim

\begin{document}

\par\q\par\q
\centerline {\bf Some relations between $t(a,b,c,d;n)$ and
$N(a,b,c,d;n)$}

$$\q$$
\centerline{Zhi-Hong Sun}
\par\q\newline
\centerline{School of Mathematical Sciences, Huaiyin Normal
University,} \centerline{Huaian, Jiangsu 223001, P.R. China}
\centerline{Email: zhihongsun@yahoo.com} \centerline{Homepage:
http://www.hytc.edu.cn/xsjl/szh}

 \abstract{Let $\Bbb Z$ and $\Bbb N$ be the set of integers
 and the set of positive integers, respectively. For
 $a,b,c,d,n\in\Bbb N$ let $N(a,b,c,d;n)$ be the number of
 representations of $n$ by $ax^2+by^2+cz^2+dw^2$, and
 let $t(a,b,c,d;n)$ be the number of
 representations of $n$ by $ax(x-1)/2+by(y-1)/2+cz(z-1)/2
 +dw(w-1)/2$ $(x,y,z,w\in\Bbb Z$). In this paper
 we reveal some connections
between $t(a,b,c,d;n)$ and $N(a,b,c,d;n)$.
 \par\q
 }
 \endabstract
\let\thefootnote\relax \footnotetext {The author
is supported by the National Natural Science Foundation of China
(grant No. 11371163).}

\section*{1. Introduction}
\par\q  Let $\Bbb Z$ and $\Bbb N$ be the set of integers
 and the set of positive integers, respectively.
 Let $\Bbb Z^4=\Bbb Z\times \Bbb Z\times \Bbb
Z\times \Bbb Z$ and $\Bbb N^4=\Bbb N\times \Bbb N\times \Bbb N\times
\Bbb N$. For  $n \in \Bbb N$ let
$$\sigma(n)=\sum_{d \mid n,d\in\Bbb N}d.$$ For convenience
 we define $\sigma(n)=0$ for $n\notin \Bbb N$. For $a,b,c,d\in\Bbb N$ and $n\in\Bbb N \cup \{0\}$ set
$$N(a,b,c,d;n)=\big|\{(x,y,z,w)\in \Bbb Z^4\ |\ n=ax^2+by^2+cz^2+dw^2
\}\big|$$ and $$t(a,b,c,d;n)=\Big|\Big\{(x,y,z,w)\in \Bbb Z^4\ |\ n\
=a\f{x(x-1)}2+ b\f{y(y-1)}2+c\f{z(z-1)}2+d\f{w(w-1)}2\Big\}\Big|.$$
The numbers $\f{x(x-1)}2\ (x\in\Bbb Z)$ are called triangular
numbers.
\par In 1828 Jacobi showed that
$$N(1,1,1,1;n)=8\sum_{d\mid n,4\nmid d}d.$$
In 1847 Eisenstein (see [D]) gave formulas for $N(1,1,1,3;n)$ and
$N(1,1,1,5;n)$. From 1859 to 1866 Liouville made about 90
conjectures on $N(a,b,c,d;n)$ in a series of papers. Most
conjectures of Liouville have been proved. See [A1, A2,
AALW1-AALW5], Cooper's survey paper [C], Dickson's historical
comments [D] and Williams' book [W3].
\par
 Let
$$t'(a,b,c,d;n)=\Big|\Big\{(x,y,z,w)\in \Bbb N^4\ |\ n=a\f{x(x-1)}2+
b\f{y(y-1)}2+c\f{z(z-1)}2+d\f{w(w-1)}2\Big\}\Big|.$$ As $\f
{x(x-1)}2=\f{(-x+1)(-x)}2$ we have
$$t(a,b,c,d;n)=16t'(a,b,c,d;n).\tag 1.1$$
In [L] Legendre stated that
$$t'(1,1,1,1;n)=\sigma(2n+1).\tag 1.2$$ In 2003,
Williams [W1] showed that
$$t'(1,1,2,2;n)=\f 14\sum_{d\mid 4n+3}
\big(d-(-1)^{\f{d-1}2}\big).\tag 1.3$$
 For $a,b,c,d\in\Bbb N$ with $5\le a+b+c+d\le 8$ let
$$C(a,b,c,d)=16+4i_1(i_1-1)i_2+8i_1i_3,$$
where $i_j$ is the number of elements in $\{a,b,c,d\}$ which are
equal to $j$. When $a+b+c+d\in\{5,6,7\}$, in 2005 Adiga, Cooper and
Han [ACH] showed that
$$C(a,b,c,d)t'(a,b,c,d;n)=N(a,b,c,d;8n+a+b+c+d).\tag 1.4$$ When
$a+b+c+d=8$, in 2008 Baruah, Cooper and Han [BCH] proved that
$$C(a,b,c,d)t'(a,b,c,d;n)=N(a,b,c,d;8n+8)-N(a,b,c,d;2n+2).\tag 1.5$$
 In 2009,
Cooper [C] determined $t'(a,b,c,d;n)$ for $(a,b,c,d)=(1,1,1,3),\
(1,3,3,3),$ $(1,2,2,3),\ (1,3,6,6),\ (1,3,4,4),\ (1,1,2,6)$ and
$(1,3,12,12)$.
\par In [WS], Wang and Sun obtained
explicit
 formulas for $t(a,b,c,d;n)$ in the cases
 $(a,b,c,d)=(1,2,2,4),\ (1,2,4,4),\ (1,1,4,4),\ (1,4,4,4)$,
 $(1,3,3,9)$,
 $(1,1,9,9),\ (1,9,9,$ $9)$, $(1,1,1,9)$, $(1,3,9,9)$ and $(1,1,3,9).$

\par Ramanujan's theta functions $\varphi(q)$ and $\psi(q)$ are defined
by
$$\varphi(q)=\sum_{n=-\infty}^{\infty}q^{n^2}=1+2\sum_{n=1}^{\infty}
q^{n^2}\qtq{and} \psi(q)=\sum_{n=0}^{\infty}q^{n(n+1)/2}\ (|q|<1).$$
It is evident that for $|q|<1$,
$$\align&\sum_{n=0}^{\infty}N(a,b,c,d;n)q^{n}=\varphi(q^a)
\varphi(q^b)\varphi(q^c)\varphi(q^d),\tag 1.6
\\&\sum_{n=0}^{\infty}t'(a,b,c,d;n)q^{n}=\psi(q^a)\psi(q^b)
\psi(q^c)\psi(q^d).\tag 1.7\endalign$$
From [BCH, Lemma 4.1] or [Be]
we know that for $|q|<1$,
 $$\align &\psi(q)^2=\varphi(q)\psi(q^2).\tag
 1.8\\&\varphi(q)=\varphi(q^4)+2q\psi(q^8),\tag 1.9
 \\&\varphi(q)^2=\phq 2^2+4q\psq 4^2,\tag 1.10
 \\&\psi(q)\psi(q^3)=\varphi(q^6)\psi(q^4)+q
 \varphi(q^2)\psi(q^{12}).\tag 1.11
 \\&\endalign$$
 By (1.9), for $k\in\Bbb N$,
$$\varphi(q^k)=\varphi(q^{4k})+2q^k\psi(q^{8k})
=\varphi(q^{16k})+2q^{4k}\psi(q^{32k})+2q^k\psi(q^{8k}).\tag 1.12$$
\par In this paper, using (1.8)-(1.12) we reveal
some connections between $t(a,b,c,d;n)$ and $N(a,b,c,d;n)$. Suppose
$k,m\in\{0,1,2,\ldots\}$, $a,n\in\Bbb N$ and $2\nmid a$. We show
that
$$t(a,b,c,d;n)=\frac
23\big(N(a,b,c,d;8n+a+b+c+d)-N(a,b,c,d;2n+(a+b+c+d)/4)\big)$$ for
$(a,b,c,d)= (a,a,2a,8m+4)$ and $(a,3a,4k+2,4m+2)$ with $k\equiv
m\pmod 2$. For $2\nmid ak$ we show that
$$t(a,3a,k,k;m)=\f 23N(a,3a,2k,2k;8m+4a+2k).$$
 For $n\e k+\f{a-1}2\mod
2$, we prove that
$$t(a,3a,8k+4,4m+2;n)=\f 23N(a,3a,8k+4,4m+2;8n+4m+8k+4a+6).$$
Let $a,k,m\in\Bbb N$ with $2\nmid a$. We also state that
$$\align &t(a,a,2a,4k;4m+3a)=4t(a,2a,4a,k;m),\
t(a,a,6a,4k;4m+3a)=2t(a,a,6a,k;m),
\\&t(a,a,8a,2k;2m)=t(a,2a,2a,k;m)\qtq{and}
t(a,a,8a,2k;2m+a)=2t(a,4a,4a,k;m).\endalign$$
 In addition, we give explicit formulas for
$t(1,3,3,6;n),\ t(1,1,8,8;n)$ and $t(1,1,4,8;n)$. We also pose many
conjectures on the relations between $t(a,b,c,d;n)$ and
$N(a,b,c,d;n)$.

\section*{2. Main results}
\pro{Lemma 2.1} Suppose $a,k,m,n\in\Bbb N$ and $2\nmid a$. Then
$$N(a,a,2k,2m;2n)=N(a,a,k,m;n).$$
\endpro
Proof. Suppose $|q|<1$.  Using (1.10) we see that
$$\align&\sum_{n=0}^{\infty}N(a,a,2k,2m;n)q^{n}
\\&=\varphi(q^{a})^2\varphi(q^{2k})\varphi(q^{2m})
=\big(\varphi(q^{2a})^2+4q^a\psi(q^{4a})^2\big)
\varphi(q^{2k})\varphi(q^{2m}).\endalign$$ Extracting the even
powers we obtain
$$\sum_{n=0}^{\infty}N(a,a,2k,2m;2n)q^{2n}
=\varphi(q^{2a})^2\varphi(q^{2k})\varphi(q^{2m}).$$
 Replacing
$q$ with $q^{1/2}$ in the above formula we obtain
$$\sum_{n=0}^{\infty}N(a,a,2k,2m;2n)q^{n}
=\varphi(q^a)^2\varphi(q^k)\varphi(q^m)
=\sum_{n=0}^{\infty}N(a,a,k,m;n)q^{n} .$$ Comparing the coefficients
of $q^n$ on both sides we obtain the result.
 \pro{Lemma
2.2} For $|q|<1$ we have
$$\varphi(q)^3=\phq 4^3+6q\phq 4\psq 4^2+12q^2\psq 4^2\psq 8+8q^3
\psq 8^3.$$
\endpro
Proof. By (1.8) and (1.9),
$$\align \varphi(q)^3&=(\phq 4+2q\psq 8)^3
\\&=\phq 4^3+6q\phq 4\psq 8(\phq 4+2q\psq 8)+8q^3\psq 8^3
\\&=\phq 4^3+6q\psq 4^2(\phq 4+2q\psq 8)+8q^3\psq 8^3.
\endalign$$
This yields the result.

\pro{Theorem 2.1} Let $a\in\{1,3,5,\ldots\}$ and
$m\in\{0,1,2,\ldots\}$. For $n\in\Bbb N$ we have
$$\align &t(a,a,2a,8m+4;n)
\\&=\f 23(N(a,a,a,4m+2;4n+4m+2a+2)-N(a,a,a,4m+2;n+m+(a+1)/2))
\\&=\f
23(N(a,a,2a,8m+4;8n+8m+4a+4)-N(a,a,2a,8m+4;2n+2m+a+1)).
\endalign$$\endpro

Proof. Suppose $|q|<1$.  By Lemma 2.2,
 $$\aligned&\sum_{n=0}^{\infty}N(a,a,a,4m+2;n)q^{n}
 \\&=\varphi(q^a)^3\varphi(q^{4m+2})
 \\&=\big(
 \phq {4a}^3+6q^a\phq {4a}\psq {4a}^2+12q^{2a}\psq {4a}^2\psq {8a}
 +8q^{3a}\psq {8a}^3\big)
\\&\q\times\big(\varphi(q^{4(4m+2)})+2q^{4m+2}
 \psi(q^{8(4m+2)}) \big)
 \endaligned\tag 2.1$$
For any $r,s\in\Bbb N$ the power series expansions of $\phq{8s}^r$
and $\psq{8s}^r$ are of the form $\sum_{n=0}^{\infty}b_nq^{8n}$.
Hence from (2.1) we deduce that
$$\align&\sum_{n=0}^{\infty}N(a,a,a,4m+2;4n)q^{4n}
\\&=\varphi(q^{4a})^3\varphi(q^{16m+8})
+12q^{2a}\psq{4a}^2\psq{8a}\cdot 2q^{4m+2}\psi(q^{32m+16}).
\endalign$$
Replacing $q$ with $q^{1/4}$ in the above formula we obtain
$$\align&\sum_{n=0}^{\infty}N(a,a,a,4m+2;4n)q^n
\\&=\phq{a}^3\phq{4m+2}+24q^{m+(a+1)/2}
\psq{a}^2\psq{2a}\psq{8m+4}
\\&=\sum_{n=0}^{\infty}N(a,a,a,4m+2;n)q^n+24q^{m+(a+1)/2}
\sum_{n=0}^{\infty}t'(a,a,2a,8m+4;n)q^n.
\endalign$$
Now comparing the coefficients of $q^{n+m+(a+1)/2}$ on both sides
and then applying (1.1) and Lemma 2.1 we obtain
$$\align&\f 32t(a,a,2a,8m+4;n)\\&=24t'(a,a,2a,8m+4;n)
\\&=N(a,a,a,4m+2;4n+4m+2a+2)-N(a,a,a,4m+2;n+m+(a+1)/2)
\\&=N(a,a,2a,8m+4;8n+8m+4a+4)-N(a,a,2a,8m+4;2n+2m+2a+2).
\endalign$$
This is the result.

\pro{Lemma 2.3} For $|q|<1$ we have
$$\align\varphi(q)\phq 3&=
\phq {16}\phq{48}+4q^{16}\psq{32}\psq{96}+2q \phq {48}\psq
8+2q^3\phq{16}\psq{24} \\&\q+6q^4\psq 8\psq{24}+4q^{13}\psq
8\psq{96}+4q^7\psq{24}\psq{32}.\endalign$$
\endpro
Proof. By (1.12),
$$\align\varphi(q)\phq 3&=\big(\phq{16}+2q^4\psq{32}+2q\psq
8\big)\big(\phq{48}+2q^{12}\psq{96}+2q^3\psq{24}\big)
\\&=\phq {16}\phq{48}+4q^{16}\psq{32}\psq{96}+2q \phq {48}\psq
8+2q^3\phq{16}\psq{24}
\\&\q+2q^4\big(\phq{48}\psq{32}+q^8\phq{16}\psq{96}+
2\psq 8\psq{24}\big) \\&\q+4q^{13}\psq
8\psq{96}+4q^7\psq{24}\psq{32}.\endalign$$
 Note that $\phq{48}\psq{32}+q^8\phq{16}\psq{96}=\psq 8\psq{24}$ by
 (1.11). We then obtain the result.

\pro{Theorem 2.2} Let $a\in\{1,3,5,\ldots\}$,
$k,m\in\{0,1,2,\ldots\}$ and $k\e m\mod 2$. For $n\in\Bbb N$ we have
$$\align &t(a,3a,4k+2,4m+2;n)\\&=\f
23(N(a,3a,4k+2,4m+2;8n+4m+4k+4a+4)\\&\q-N(a,3a,4k+2,4m+2;2n+m+k+a+1)).
\endalign$$\endpro

Proof. Suppose $|q|<1$. Using Lemma 2.3 and (1.12) we see that
 $$\align &\sum_{n=0}^{\infty}N(a,3a,4k+2,4m+2;n)q^{n}
\\&=\varphi(q^a)\varphi(q^{3a})\varphi(q^{4k+2})
\varphi(q^{4m+2})
\\&=\varphi(q^a)\varphi(q^{3a})\big(\varphi(q^{4(4k+2)})+2q^{4k+2}
 \psi(q^{8(4k+2)})\big)\big(\varphi(q^{4(4m+2)})+2q^{4m+2}
 \psi(q^{8(4m+2)})\big)
\\&=\big(\phq {16a}\phq{48a}+4q^{16a}\psq{32a}\psq{96a}
+2q^a \phq {48a}\psq {8a}+2q^{3a}\phq{16a}\psq{24a}
\\&\q+6q^{4a}\psq {8a}\psq{24a}+4q^{13a}\psq
{8a}\psq{96a}+4q^{7a}\psq{24a}\psq{32a}\big)
\\&\q\times \big(\varphi(q^{16k+8})\varphi(q^{16m+8})+
2q^{4m+2}\varphi(q^{16k+8})\psi(q^{32m+16})\\&\q+2q^{4k+2}
\psi(q^{32k+16})\varphi(q^{16m+8})
+4q^{4k+4m+4}\psi(q^{32k+16})\psi(q^{32m+16})\big).
\endalign$$
For any $r,s\in\Bbb N$ the power series expansions of $\phq{8s}^r$
and $\psq{8s}^r$ are of the form $\sum_{n=0}^{\infty}b_nq^{8n}$.
Thus,
 from the above we deduce that
$$\align&\sum_{n=0}^{\infty}N(a,3a,4k+2,4m+2;8n)
q^{8n}
\\&=\big(\phq{16a}\phq{48a}
+4q^{16a}\psq{32a}\psq{96a}\big)\phq{16k+8}\phq{16m+8}
\\&\q+24q^{4(k+m+a+1)}\psq{8a}\psq{24a}\psq{32k+16}\psq{32m+16}.
\endalign$$
 Replacing $q$ with
$q^{1/8}$ in the above formula we obtain
$$\align&\sum_{n=0}^{\infty}N(a,3a,4k+2,4m+2;8n)
q^n
\\&=(\phq{2a}\phq{6a}
+4q^{2a}\psq{4a}\psq{12a})\phq{2k+1}\phq{2m+1}
\\&\q+24q^{(k+m+a+1)/2}\psq{a}\psq{3a}\psq{4k+2}\psq{4m+2}.
\endalign$$
On the other hand, using (1.9) we see that
$$\align&\sum_{n=0}^{\infty}N(a,3a,4k+2,4m+2;n)q^{n}
\\&=\varphi(q^{a})\varphi(q^{3a})\varphi(q^{4k+2})
\varphi(q^{4m+2})
\\&=\big(\varphi(q^{4a})+2q^a\psi(q^{8a})\big)
\big(\varphi(q^{12a})+2q^{3a}\psi(q^{24a})\big)
 \varphi(q^{4k+2})
\varphi(q^{4m+2}).
\endalign$$
Extracting the even powers we obtain
$$\align&\sum_{n=0}^{\infty}N(a,3a,4k+2,4m+2;2n)q^{2n}
\\&=\big(\varphi(q^{4a})\varphi(q^{12a})+
4q^{4a}\psi(q^{8a})\psi(q^{24a})\big) \varphi(q^{4k+2})
\varphi(q^{4m+2}).
\endalign$$
Replacing $q$ with $q^{1/2}$ we then obtain
$$\align&\sum_{n=0}^{\infty}N(a,3a,4k+2,4m+2;2n)q^{n}
\\&=\big(\varphi(q^{2a})\varphi(q^{6a})+
4q^{2a}\psi(q^{4a})\psi(q^{12a})\big) \varphi(q^{2k+1})
\varphi(q^{2m+1}).
\endalign$$
Hence
$$\align&\sum_{n=0}^{\infty}(N(a,3a,4k+2,4m+2;8n)
-N(a,3a,4k+2,4m+2;2n))q^{n}
\\&=24q^{(k+m+a+1)/2}\psi(q^a)\psi(q^{3a})
\psi(q^{4k+2})\psi(q^{4m+2}).
\\&=24q^{(k+m+a+1)/2}\sum_{n=0}^{\infty}t'(a,3a,4k+2,4m+2;n)q^n
\\&=\f 32q^{(k+m+a+1)/2}
\sum_{n=0}^{\infty}t(a,3a,4k+2,4m+2;n)q^n.\endalign$$
 Comparing the coefficients of
 $q^{n+(k+m+a+1)/2}$ yields the result.

\pro{Theorem 2.3} Let $a,k\in\Bbb N$ with $2\nmid ak$. For $m\in\Bbb
N$ we have
$$t(a,3a,k,k;m)=\f 23N(a,3a,2k,2k;8m+4a+2k).$$
\endpro
Proof. By (1.10),
$$\phq{2k}^2=\phq{4k}^2+4q^{2k}\psq{8k}^2=\phq{8k}^2
+4q^{4k}\psq{16k}^2+4q^{2k}\psq{8k}^2.$$ Thus, applying Lemma 2.3 we
see that
$$\align&\sum_{n=0}^{\infty}N(a,3a,2k,2k;n)q^n
\\&=\phq a\phq{3a}\phq{2k}^2
\\&=\big(\phq {16a}\phq{48a}+4q^{16a}\psq{32a}\psq{96a}
+2q^a \phq {48a}\psq {8a}+2q^{3a}\phq{16a}\psq{24a}
\\&\q+6q^{4a}\psq {8a}\psq{24a}+4q^{13a}\psq
{8a}\psq{96a}+4q^{7a}\psq{24a}\psq{32a}\big) \\&\q\times
\big(\phq{8k}^2 +4q^{4k}\psq{16k}^2+4q^{2k}\psq{8k}^2\big).
\endalign$$
For any $r,s\in\Bbb N$ the power series expansions of
$\varphi(q^s)^r$ and $\psi(q^s)^r$ are of the form
$\sum_{n=0}^{\infty}b_nq^{8n}$. Thus, from the above we deduce that
$$\align &\sum_{m=0}^{\infty}N(a,3a,2k,2k;8m+4a+2k)q^{8m+4a+2k}
\\&=6q^{4a}\psq{8a}\psq{24a}\cdot 4q^{2k}\psq{8k}^2
=24q^{4a+2k}\psq{8a}\psq{24a}\psq{8k}^2.\endalign$$ Replacing $q$
with $q^{1/8}$ we get
$$\align &\sum_{m=0}^{\infty}N(a,3a,2k,2k;8m+4a+2k)q^m
\\&=24\psq a\psq{3a}\psq k^2
=24\sum_{m=0}^{\infty}t'(a,3a,k,k;m)q^m =\f
{24}{16}\sum_{m=0}^{\infty}t(a,3a,k,k;m)q^m.
\endalign$$
Comparing the coefficients of $q^m$ on both sides yields the result.
\pro{Corollary 2.1} For $n\in\Bbb N$ we have
$$N(2,2,3,9;8n+6)=\f 35N(1,1,3,9;8n+6).$$
\endpro
Proof. Taking $a=3$ and $k=1$ in Theorem 2.3 we see that
$t(1,1,3,9;m)=\f 23N(2,2,3,9;8m+14).$ On the other hand, from [WS,
the proof of Theorem 2.3] we know that $t(1,1,3,9;m)=\f
25N(1,1,3,9;8m+14)$. Thus, the result follows.

\pro{Theorem 2.4} Let $a\in\{1,3,5,\ldots\}$,
$k,m\in\{0,1,2,\ldots\}$ and $n\in\Bbb N$. If $n\e k+\f{a-1}2\mod
2$, then
$$t(a,3a,8k+4,4m+2;n)=\f 23N(a,3a,8k+4,4m+2;8n+4m+8k+4a+6).$$
\endpro
Proof. Suppose $|q|<1$. Using Lemma 2.3 and (1.12) we see that
 $$\align &\sum_{n=0}^{\infty}N(a,3a,8k+4,4m+2;n)q^{n}
\\&=\varphi(q^{a})\varphi(q^{3a})\varphi(q^{8k+4})\varphi(q^{4m+2})
\\&=\varphi(q^{a})\varphi(q^{3a})\big(\varphi(q^{4(8k+4)})+2q^{8k+4}
 \psi(q^{8(8k+4)})\big)\big(\varphi(q^{4(4m+2)})+2q^{4m+2}
 \psi(q^{8(4m+2)})\big)
\\&=\big(\phq {16a}\phq{48a}+4q^{16a}\psq{32a}\psq{96a}
+2q^a \phq {48a}\psq {8a}+2q^{3a}\phq{16a}\psq{24a}
\\&\q+6q^{4a}\psq {8a}\psq{24a}+4q^{13a}\psq
{8a}\psq{96a}+4q^{7a}\psq{24a}\psq{32a}\big)
\\&\q\times \big(\varphi(q^{32k+16})\varphi(q^{16m+8})+
2q^{4m+2}\varphi(q^{32k+16})\psi(q^{32m+16})\\&\q+2q^{8k+4}
\psi(q^{64k+32})\varphi(q^{16m+8})
+4q^{8k+4m+6}\psi(q^{64k+32})\psi(q^{32m+16})\big).
\endalign$$
 For any $r,s\in\Bbb N$ the power series expansions of $\phq{8s}^r$
and $\psq{8s}^r$ are of the form $\sum_{n=0}^{\infty}b_nq^{8n}$.
Thus, from the above and the fact that $4m+2\e 4-2(-1)^m\mod 8$ we
deduce that
$$\align&\sum_{n=0}^{\infty}N(a,3a,8k+4,4m+2;8n+4-2(-1)^m)
q^{8n+4-2(-1)^m}
\\&=\big(\phq {16a}\phq{48a}+4q^{16a}\psq{32a}\psq{96a}\big)
\cdot2q^{4m+2}\varphi(q^{32k+16})\psi(q^{32m+16})
\\&\q+6q^{4a}\psq {8a}\psq{24a}\cdot
4q^{8k+4m+6}\psi(q^{64k+32})\psi(q^{32m+16})
\endalign$$ and so
$$\align&\sum_{n=0}^{\infty}N(a,3a,8k+4,4m+2;8n+4-2(-1)^m)q^{8n}
\\&=2q^{8[m/2]}
\varphi(q^{32k+16})\psi(q^{32m+16})\big(\phq
{16a}\phq{48a}+4q^{16a}\psq{32a}\psq{96a}\big)
\\&\q+24q^{8(k+[m/2]+(a+1)/2)}
\psq {8a}\psq{24a} \psi(q^{64k+32})\psi(q^{32m+16}),
\endalign$$
where $[x]$ is the greatest integer not exceeding $x$.
 Replacing $q$ with
$q^{1/8}$ in the above formula we obtain
$$\align&\sum_{n=0}^{\infty}N(a,3a,8k+4,4m+2;8n+4-2(-1)^m)q^n
\\&=2q^{[m/2]}
\varphi(q^{4k+2})\psi(q^{4m+2})\big(\phq
{2a}\phq{6a}+4q^{2a}\psq{4a}\psq{12a}\big)
\\&\q+24q^{k+[m/2]+(a+1)/2}
\psq {a}\psq{3a} \psi(q^{8k+4})\psi(q^{4m+2}),
\endalign$$
 Suppose that $n\e k+(a-1)/2\mod 2$. Then $n+k+[m/2]+(a+1)/2
 \e [m/2]+1\mod 2$.
 Now comparing the coefficients of
 $q^{n+k+[m/2]+(a+1)/2}$ in
the above expansion we obtain
$$\align&N(a,3a,8k+4,4m+2;8(n+k+[m/2]+(a+1)/2)+4-2(-1)^m)
\\&=24t'(a,3a,8k+4,
4m+2;n)=\f32t(a,3a,8k+4,4m+2;n).\endalign$$ This yields the result.

\pro{Theorem 2.5} Let $a,k\in\Bbb N$ with $2\nmid a$. For $m\in\Bbb
N$ we have
$$t(a,a,6a,4k;4m+3a)=2t(a,a,6a,k;m).$$
\endpro
Proof. Suppose $|q|<1$. Using (1.8)-(1.12) we see that
$$\aligned &\sum_{n=0}^{\infty}t'(a,a,6a,4k;n)q^n
\\&=\psq a^2\psq {6a}\psq{4k}
=\phq a\psq{2a}\psq{6a}\psq{4k}
\\&=\big(\phq{4a}+2q^a\psq{8a}\big)\big(\phq{12a}\psq{8a}+q^{2a}\phq{4a}\psq{24a}
\big)\psq{4k}
\\&=\big(\phq{4a}\phq{12a}\psq{8a}+2q^a\psq{8a}^2\phq{12a}
+q^{2a}\phq{4a}^2\psq{24a}
\\&\q+2q^{3a}\phq{4a}
\psq{8a}\psq{24a}\big)\psq{4k} .\endaligned\tag 2.2$$
 Extracting the
powers of $q^{4m+3a}$ we get
$$\align &\sum_{m=0}^{\infty}t'(a,a,6a,4k;4m+3a)q^{4m+3a}
\\&=2q^{3a}\phq{4a} \psq{8a}\psq{24a}\psq{4k} =2q^{3a}
\psq{4a}^2\psq{24a}\psq{4k}.\endalign$$ Replacing $q$ with $q^{1/4}$
we deduce that
$$\sum_{m=0}^{\infty}t'(a,a,6a,4k;4m+3a)q^m
=2\psq{a}^2\psq{6a}\psq{k}=2\sum_{m=0}^{\infty} t'(a,a,6a,k;m)q^m.$$
Hence
$$\align t(a,a,6a,4k;4m+3a)&=16t'(a,a,6a,4k;4m+3a)=32t'(a,a,6a,k;m)
\\&=2t(a,a,6a,k;m)\endalign$$ as asserted.

\pro{Theorem 2.6} Let $a,k\in\Bbb N$ with $2\nmid a$. For $m\in\Bbb
N$ we have
$$t(a,a,2a,4k;4m+3a)=4t(a,2a,4a,k;m).$$
\endpro
Proof. Suppose $|q|<1$. Using (1.8) and (1.12) we see that
$$\align &\sum_{n=0}^{\infty}t'(a,a,2a,4k;n)q^n
\\&=\psq a^2\psq {2a}\psq{4k}
=\phq a\psq{2a}^2\psq{4k}
\\&=\phq a\phq{2a}\psq{4a}\psq{4k}
\\&=
\big(\phq{4a}+2q^{a}\psq{8a}\big)\big(\phq{8a}+2q^{2a}
\psq{16a}\big)\psq{4a}\psq{4k}
\\&=\big(\phq{4a}\phq{8a}+2q^a\phq{8a}\psq{8a}
+2q^{2a}\phq{4a}\psq{16a}
\\&\q+4q^{3a}\psq{8a}
\psq{16a}\big)\psq{4a}\psq{4k} .\endalign$$
 Extracting the
powers of $q^{4m+3a}$ we get
$$\sum_{m=0}^{\infty}t'(a,a,2a,4k;4m+3a)q^{4m+3a}
=4q^{3a}\psq{8a} \psq{16a}\psq{4a}\psq{4k} .$$ Replacing $q$ with
$q^{1/4}$ we see that
$$\align &\sum_{m=0}^{\infty}t'(a,a,2a,4k;4m+3a)q^m
\\&=4\psq{a}\psq {2a}\psq {4a}\psq k=4\sum_{m=0}^{\infty}
t'(a,2a,4a,k;m)q^m.\endalign$$ Hence
$$\align t(a,a,2a,4k;4m+3a)&=16t'(a,a,2a,4k;4m+3a)
=64t'(a,2a,4a,k;m)\\&=4t(a,2a,4a,k;m), \endalign$$
 which completes the proof.

 \pro{Theorem 2.7} Let $a,k\in\Bbb N$ with $2\nmid a$. For $n\in\Bbb
N$ we have
$$t(a,a,8a,2k;2n)=t(a,2a,2a,k;n)$$
and
$$t(a,a,8a,2k;2n+a)=2t(a,4a,4a,k;n).$$
\endpro
Proof. Suppose $|q|<1$. Using (1.8) and (1.9) we see that
$$\aligned \sum_{n=0}^{\infty}t'(a,a,8a,2k;n)q^n
&=\psq a^2\psq {8a}\psq{2k} =\phq a\psq{2a}\psq {8a}\psq{2k}
\\&=\psq{2a}(\phq{4a}+2q^a\psq{8a})\psq {8a}\psq{2k}.
\endaligned\tag 2.3$$
Extracting the even powers we see that
$$\sum_{n=0}^{\infty}t'(a,a,8a,2k;2n)q^{2n}
=\psq{2a}\phq{4a}\psq {8a}\psq{2k}=\psq{2a}\psq{4a}^2\psq{2k}.
$$
Replacing $q$ with $q^{1/2}$ we then get
$$\sum_{n=0}^{\infty}t'(a,a,8a,2k;2n)q^n
=\psq{a}\psq{2a}^2\psq{k} =\sum_{n=0}^{\infty}t'(a,2a,2a,k;n)q^n.
$$
Hence
$$t(a,a,8a,2k;2n)=16t'(a,a,8a,2k;2n)=16t'(a,2a,2a,k;n)
=t(a,2a,2a,k;n).$$ On the other hand, extracting the odd powers in
(2.3) we see that
$$\sum_{n=0}^{\infty}t'(a,a,8a,2k;2n+a)q^{2n+a}
=\psq{2a}\cdot 2q^a\psq{8a}\psq {8a}\psq{2k}.$$ Replacing $q$ with
$q^{1/2}$ we then get
$$\sum_{n=0}^{\infty}t'(a,a,8a,2k;2n+a)q^n
=2\psq{a}\psq{4a}^2\psq{k} =\sum_{n=0}^{\infty}2t'(a,4a,4a,k;n)q^n.
$$
Hence
$$t(a,a,8a,2k;2n+a)=16t'(a,a,8a,2k;2n+a)=32t'(a,4a,4a,k;n)
=2t(a,4a,4a,k;n).$$ We are done.

\pro{Theorem 2.8} For $n\in\Bbb N$ we have
$$t(1,1,8,8;n)=\sigma(4n+9)-(2-(-1)^n)\sum\Sb (x,y)\in\Bbb Z
\times\Bbb Z, x\e 1\pmod 4
\\4n+9=x^2+4y^2\endSb x.$$
\endpro
Proof. By Theorem 2.7, $t(1,1,8,8;n)=t(1,2,2,4;n/2)$ for even $n$,
and $t(1,1,8,8;n)=2t(1,4,4,4;(n-1)/2)$ for odd $n$. Now applying
[WS, Theorems 3.2 and 3.3] yields the result.

\pro{Theorem 2.9} For $n\in\Bbb N$ we have
$$t(1,1,4,8;n)=2(-1)^n\sum_{d\mid 4n+7}d\Ls 2d
-(1-(-1)^n)\sum\Sb (x,y)\in\Bbb Z\times \Bbb Z,x\e 1\mod 4
\\4n+7=x^2+2y^2\endSb x.$$
\endpro
Proof. By Theorem 2.7, $t(1,1,4,8;n)=t(1,2,2,2;n/2)$ for even $n$,
and $t(1,1,4,8;n)=2t(1,2,4,4;(n-1)/2)$ for odd $n$. It is easily
seen that
$$\align t(1,2,2,2;m)&=\big|\big\{(x,y,z,w)\in \Bbb Z^4\ \big|\
8m+7=x^2+2y^2+2z^2+2w^2, \ 2\nmid xyzw\big\}\big|
\\&=\big|\big\{(x,y,z,w)\in \Bbb Z^4\ \big|\
8m+7=x^2+2y^2+2z^2+2w^2\big\}\big| \\&=N(1,2,2,2;m).\endalign$$
 Now applying [W2, (1.4)] and [WS, Theorem 3.4] yields the
result.

 \pro{Lemma 2.4} For $|q|<1$ we have
$$\varphi(q)^2=\phq 8^2+4q^4\psq{16}^2+4q^2\psq 8^2+4q\phq{16}\psq
8+8q^5\psq 8\psq{32}.$$
\endpro
Proof. By (1.10),
$$\align \varphi(q)^2-\phq 8^2&=(\varphi(q)^2-\phq 2^2)+(\phq 2^2-\phq 4^2)+
(\phq 4^2-\phq 8^2)
\\&=4q\psq 4^2+4q^2\psq 8^2+4q^4\psq{16}^2.\endalign$$
By (1.8) and (1.9),
$$\psq 4^2=\phq 4\psq 8=(\phq{16}+2q^4\psq{32})\psq 8.$$
Now combining the above we deduce the result.

 \pro{Theorem 2.10} If $n\in\Bbb N$ and $8n+13=3^{\beta}n_1$ with
$n_1\in\Bbb N$ and $3\nmid n_1$, then
$$\align t(1,3,3,6;n)&=\f 25N(1,3,3,6;8n+13)
\\&=\f 23\Big(3^{\beta}+\Ls{n_1}3\Big)
\prod_{p\mid n_1}\f{p^{\t{\rm ord}_pn_1+1}-\sls 6p^{\t{\rm
ord}_pn_1+1}}{p-\sls 6p},
\endalign$$
where $p$ runs all distinct prime divisors of $n_1$, $\sls ap$ is
the Legendre symbol and $\t{ord}_pn_1$ is the unique nonnegative
integer $r$ such that $p^r\mid n_1$ but $p^{r+1}\nmid n_1$.
\endpro
Proof. Suppose $|q|<1$. By (1.12) and Lemma 2.4,
$$\align &\sum_{n=0}^{\infty}N(1,3,3,6;n)q^n
\\&=\varphi(q)\varphi(q^6)\varphi(q^3)^2
\\&=(\phq{16}+2q^4\psq{32}+2q\psq{8}) (\phq{24}+2q^6\psq{48})
\\&\q\times\big(\phq {24}^2+4q^{12}\psq{48}^2+4q^6\psq {24}^2
+4q^3\phq{48}\psq {24}+8q^{15}\psq {24}\psq{96}\big)
\\&=\big(\phq{16}\phq{24}+2q^4\psq{32}\phq{24}+2q\psq{8}\phq{24}
\\&\q+2q^6\phq{16}\psq{48}+4q^{10}\psq{32}\psq{48}+4q^7\psq{8}\psq{48}
\big)
\\&\q\times\big(\phq {24}^2+4q^{12}\psq{48}^2+4q^6\psq {24}^2
+4q^3\phq{48}\psq {24}+8q^{15}\psq {24}\psq{96}\big) .
\endalign$$
 For any $r,s\in\Bbb N$ the power series expansions of
$\varphi(q^s)^r$ and $\psi(q^s)^r$ are of the form
$\sum_{n=0}^{\infty}b_nq^{8n}$. Thus, from the above we deduce that
$$\align &\sum_{n=0}^{\infty}N(1,3,3,6;8n+13)q^{8n+13}
\\&=\sum_{n=0}^{\infty}N(1,3,3,6;8n+5)q^{8n+5}
\\&=2q\psq8\phq{24}\cdot 4q^{12}\psq{48}^2
+2q^6\phq{16}\psq{48}\cdot 8q^{15}\psq{24}\psq{96}
\\&\q+4q^{10}\psq{32}\psq{48}\cdot 4q^3\psq{24}\phq{48}
+4q^7\psq 8\psq{48}\cdot 4q^6\psq{24}^2
\endalign$$
and so
$$\align &\sum_{n=0}^{\infty}N(1,3,3,6;8n+13)q^{8n}
\\&=8\psq 8\phq{24}\phq{48}\psq{96}+16q^8\phq{16}\psq{48}\psq{24}
\psq{96}
\\&\q+16\psq{24}\psq{32}\psq{48}\phq{48}+16\psq 8\psq{48}\psq{24}^2.
\endalign$$
By (1.8),
$$\phq{24}\phq{48}\psq{96}=\f{\psq{24}^2}{\psq{48}}\cdot
\f{\psq{48}^2}{\psq{96}}\cdot\psq{96}=\psq{24}^2\psq{48}.$$ By
(1.11),
$$\psq{32}\phq{48}+q^8\phq{16}\psq{96}=\psq 8\psq{24}.$$
Hence,
$$\align &\sum_{n=0}^{\infty}N(1,3,3,6;8n+13)q^{8n}
\\&=8\psq 8\psq{24}^2\psq{48}+16\psq{24}
\psq{48}\psq 8\psq{24} +16\psq 8\psq{48}\psq{24}^2
\\&=40\psq 8\psq{24}^2\psq{48} .
\endalign$$
Replacing $q$ with $q^{1/8}$ we see that
$$\align &\sum_{n=0}^{\infty}N(1,3,3,6;8n+13)q^n
\\&=40\psi(q)\psq{3}^2\psq{6}=40\sum_{n=0}^{\infty}t'(1,3,3,6;n)q^n
=\f 52 \sum_{n=0}^{\infty}t(1,3,3,6;n)q^n.
\endalign$$
Hence
$$t(1,3,3,6;n)=\f 25N(1,3,3,6;8n+13).$$
Now applying [AW, Theorem 4.1] yields the remaining part.
\par\q
\par{\bf Remark 2.1} Using Maple we find
$$t(a,b,c,d;n)=\f 25N(a,b,c,d;8n+a+b+c+d)$$ for
$(a,b,c,d)=(1,1,1,2),(1,1,1,3),(1,1,2,3),(1,1,3,9),(1,3,3,3),(1,3,3,6)$
and $(1,3,9,9)$. The cases $(1,1,3,9)$ and $(1,3,9,9)$ were solved
by Wang and Sun in [WS].

 \pro{Theorem 2.11} If $n\in\Bbb N$ and
$n\e 1,2\mod 4$, then
$$t(1,1,4,6;n)=2N(1,1,4,6;2n+3).$$
\endpro
Proof. By (1.12) and Lemma 2.4,
$$\aligned &\sum_{n=0}^{\infty}N(1,1,4,6;n)q^n
\\&=\varphi(q)^2\phq 4\phq 6
\\&=\varphi(q)^2(\phq{16}+2q^4\psq{32})(\phq{24}+2q^6\psq{48})
\\&=\big(\phq 8^2+4q^4\psq{16}^2+4q^2\psq 8^2+4q\phq{16}\psq
8+8q^5\psq
8\psq{32}\big)\\&\q\times\big(\phq{16}\phq{24}+4q^{10}\psq{32}\psq{48}
+2q^4\phq{24}\psq{32}+2q^6\phq{16}\psq{48}\big).
\endaligned\tag 2.4$$
Thus,
$$\align&\sum_{m=0}^{\infty}N(1,1,4,6;8m+5)q^{8m+5}
=16q^5\psq 8\psq{32}\phq{16}\phq{24},
\\&\sum_{m=0}^{\infty}N(1,1,4,6;8m+7)q^{8m+7} =8q^7\psq
8\psq{48}\big(\phq{16}^2+4q^8\psq{32}^2\big).
\endalign$$
 From the above and (1.10) we deduce that
$$\aligned&\sum_{m=0}^{\infty}N(1,1,4,6;8m+5)q^m
=16\psi(q)\psq{4}\phq{2}\phq{3}=16\psi (q)\psq 2^2\phq 3,
\\&\sum_{m=0}^{\infty}N(1,1,4,6;8m+7)q^m =8\psi(q)
\psq{6}\big(\phq{2}^2+4q\psq{4}^2\big)=8\varphi(q)^2\psi(q)\psq{6}.
\endaligned\tag 2.5$$
On the other hand, from (2.2) we know that

$$\align&\sum_{m=0}^{\infty}t'(1,1,4,6;4m+1)q^{4m+1}
=2q\psq 4\psq 8^2\phq{12},
\\&\sum_{m=0}^{\infty}t'(1,1,4,6;4m+2)q^{4m+2}
=q^2\phq 4^2\psq 4\psq {24}.\endalign$$ Therefore,
$$\align&\sum_{m=0}^{\infty}t'(1,1,4,6;4m+1)q^m
=2\psi(q)\psq 2^2\phq{3},
\\&\sum_{m=0}^{\infty}t'(1,1,4,6;4m+2)q^m
=\varphi(q)^2\psi(q)\psq 6.\endalign$$ This together with (2.5)
yields $N(1,1,4,6;8m+5)=8t'(1,1,4,6;4m+1)$ and
$N(1,1,4,6;8m+7)=8t'(1,1,4,6;4m+2)$. To complete the proof, we
recall that $t(a,b,c,d;n)=16t'(a,b,c,d;n)$.

\pro{Theorem 2.12} For $n\in\Bbb N$
with $n\e 1\mod 4$ we have
$$t(2,2,3,9;n)=\f 43N(2,2,3,9;2n+4).$$
\endpro
Proof. Using (1.8)-(1.12) we see that for $|q|<1$,
$$\align &\sum_{n=0}^{\infty}t'(2,2,3,9;n)q^n
\\&=\psq 2^2\psq 3\psq 9
=\phq 2\psq 4(\phq {18}\psq{12}+q^3\phq 6\psq{36})
\\&=(\phq 8+2q^2\psq{16})\psq 4\big((\phq
{72}+2q^{18}\psq{144})\psq{12}
\\&\q+q^3(\phq{24}+2q^6\psq{48})\psq{36}\big).
\endalign$$
Hence
$$\align &\sum_{m=0}^{\infty}t'(2,2,3,9;4m+5)q^{4m+5}
\\&=(2q^5\phq{24}\psq{16}+2q^9\phq 8\psq{48})\psq 4\psq{36}
\\&=2q^5\psq 4\psq{12}\psq 4\psq{36}\endalign$$
and so
$$\sum_{m=0}^{\infty}t'(2,2,3,9;4m+5)q^m
=2\psi(q)^2\psq 3\psq 9=2\sum_{m=0}^{\infty}t'(1,1,3,9;m)q^m.
$$
Therefore,
$$t(2,2,3,9;4m+5)=16t'(2,2,3,9;4m+5)=32t'(1,1,3,9;m)
=2t(1,1,3,9;m).$$ Taking $a=3$ and $k=1$ in Theorem 2.3 we see that
$t(1,1,3,9;m)=\f 23N(2,2,3,9;8m+6).$ Thus,
$t(2,2,3,9;4m+5)=2t(1,1,3,9;m)=\f 43N(2,2,3,9;8m+14).$ This yields
the result.
 \pro{Theorem 2.13} For $n\in\Bbb N$ we have
$$t(1,2,2,6;n)=\f 12N(1,1,4,6;8n+11)$$
and
$$t(1,1,8,12;2n)=\f 12N(1,1,8,12;16n+22).$$
\endpro
Proof. From (2.3) and (1.8) we know that
$$\align&\sum_{m=0}^{\infty}N(1,1,4,6;8m+3)q^{8m+3}
\\&=32q^{11}\phq{16}\psq 8\psq{32}\psq{48} =32q^{11}\psq
8\psq{16}^2\psq{48}.\endalign$$ Hence
$$\align&\sum_{m=0}^{\infty}N(1,1,4,6;8m+3)q^m
\\&=32q\psi(q)\psq 2^2\psq
6=32\sum_{n=0}^{\infty}t'(1,2,2,6;n)q^{n+1}
=2\sum_{n=0}^{\infty}t(1,2,2,6;n)q^{n+1}.
\endalign$$
Now comparing the coefficients of $q^{n+1}$ we obtain
$N(1,1,4,6;8n+11)=2t(1,2,2,6;n)$. Applying Lemma 2.1 and Theorem 2.7
we see that
$$t(1,1,8,12;2n)=t(1,2,2,6;n)=\f 12N(1,1,4,6;8n+11)
=\f 12N(1,1,8,12;16n+22).$$ We are done.

\pro{Theorem 2.14} For $m\in\Bbb N$ we have
$$t(1,1,6,24;4m+1)=2t(2,2,3,3;m)$$
and
$$t(1,1,6,24;4m)=t(1,1,3,3;m)=2^{\alpha+4}\sigma(m_1),$$
where $\alpha$ and $m_1$  are given by $m+1=2^{\alpha}3^{\beta}m_1$
and $\t{gcd}(m_1,6)=1$.
\endpro
Proof. Taking $a=1$ and $k=6$ in (2.2) we see that
$$\align &\sum_{n=0}^{\infty}t'(1,1,6,24;n)q^n
\\&=\big(\phq{4}\phq{12}\psq{8}+2q\psq{8}^2\phq{12}
+q^{2}\phq{4}^2\psq{24}
\\&\q+2q^{3}\phq{4}
\psq{8}\psq{24}\big)\psq{24} .\endalign$$ Hence, using (1.8) we see
that
$$\sum_{m=0}^{\infty}t'(1,1,6,24;4m)q^{4m}
=\phq{4}\psq{8}\phq{12}\psq{24}=\psq 4^2\psq{12}^2$$ and
$$\sum_{m=0}^{\infty}t'(1,1,6,24;4m+1)q^{4m+1}
=2q\psq{8}^2\phq{12}\psq{24}=2q\psq 8^2\psq{12}^2.$$ Therefore,
$$\align &\sum_{m=0}^{\infty}t'(1,1,6,24;4m)q^m
=\psi(q)^2\psq{3}^2=\sum_{m=0}^{\infty}t'(1,1,3,3;m)q^m,
\\&\sum_{m=0}^{\infty}t'(1,1,6,24;4m+1)q^m
=2\psq 2^2\psq{3}^2=2\sum_{m=0}^{\infty}t'(2,2,3,3;m)q^m.
\endalign$$
Since $t(a,b,c,d;n)=16t'(a,b,c,d;n)$, we obtain
$t(1,1,6,24;4m+1)=2t(2,2,3,3;m)$ and $t(1,1,6,24;4m)=t(1,1,3,3;m)$.
Now applying [WS, Lemma 4.1] we deduce the theorem.
\par\q
\par In conclusion, we pose the following conjectures.
\pro{Conjecture 2.1} Let $n\in\Bbb N$ with $n\e 0,3\mod 4$. Then
$$t(1,1,4,6;n)=\f 23N(1,1,4,6;8n+12)-N(1,1,4,6;2n+3).$$
\endpro

\pro{Conjecture 2.2} Let $n\in\Bbb N$. If $3\mid n$, then
$$t(1,1,8,12;n)=\f 12N(1,1,8,12;8n+22).$$
\endpro

\pro{Conjecture 2.3} Let $m\in\Bbb N$. Then
$$t(1,3,8,8;3m)=\f 13N(1,3,8,8;24m+20)-2N(1,3,8,8;6m+5).$$
\endpro

\pro{Conjecture 2.4} Let $n\in\Bbb N$ with $n\e 0\mod 6$. Then
$$t(1,2,3,8;n)=\f 23N(1,2,3,8;8n+14)-2N(1,2,3,8;4n+7).$$
\endpro

\pro{Conjecture 2.5} Let $n\in\Bbb N$ with $n\e 0,2\mod 8$. Then
$$t(1,2,4,17;n)=4N(1,2,4,17;n+3).$$
\endpro

\pro{Conjecture 2.6} Let $n\in\Bbb N$. If $n\e 2,3\mod 5$, then
$$t(1,1,5,8;n)=\f 12N(1,1,5,8;8n+15).$$
\endpro

\pro{Conjecture 2.7} Let $n\in\Bbb N$. If $n\e 0,3,4,6,7\mod 9$,
then
$$t(1,1,8,9;n)=\f 12N(1,1,8,9;8n+19).$$
\endpro

\pro{Conjecture 2.8} Let $n\in\Bbb N$. If $n\e 0,4,7,8,9,10\mod
{13}$, then
$$t(1,1,8,13;n)=\f 12N(1,1,8,13;8n+23).$$
\endpro

\pro{Conjecture 2.9} Let $n\in\Bbb N$. If $n\e 0,3,5,6,7\mod {11}$,
then
$$t(1,1,4,11;n)=\f 13N(1,1,4,11;8n+17).$$
\endpro

\pro{Conjecture 2.10} Let $n\in\Bbb N$. If $n\e 0,1,2,4,7\mod {11}$,
then
$$t(1,1,2,22;n)=\f 13N(1,1,2,22;8n+26).$$
\endpro

\pro{Conjecture 2.11} Let $n\in\Bbb N$ with $n\e 1\mod 3$. Then
$$t(1,3,12,36;n)=\f 12N(1,3,12,36;8n+52)-2N(1,3,12,36;2n+13).$$
\endpro

\pro{Conjecture 2.12} Let $n\in\Bbb N$ with $n\e 1\mod 4$. Then
$$t(3,5,20,32;n)=\f 12N(3,5,20,32;8n+60)-2N(3,5,20,32;2n+15).$$
\endpro

\pro{Conjecture 2.13} Let $n\in\Bbb N$ with $n\e 1\mod 4$. Then
$$t(1,6,15,18;n)=\f 23N(1,6,15,18;8n+40)-2N(1,6,15,18;2n+10).$$
\endpro

\pro{Conjecture 2.14} Let $n\in\Bbb N$ with $n\e 1\mod 3$. Then
$$t(1,6,18,27;n)=\f 23N(1,6,18,27;8n+52)-2N(1,6,18,27;2n+13).$$
\endpro

\pro{Conjecture 2.15} Let $n\in\Bbb N$ with $n\e 1\mod 3$. Then
$$t(1,8,9,18;n)=\f 23N(1,8,9,18;8n+36)-2N(1,8,9,18;2n+9).$$
\endpro

\pro{Conjecture 2.16} Let $n\in\Bbb N$ with $4\mid n$. Then
$$t(1,7,10,30;n)=\f 23N(1,7,10,30;8n+48)-2N(1,7,10,30;2n+12).$$
\endpro

\pro{Conjecture 2.17} Let $n\in\Bbb N$ with $n\e 3\mod 4$. Then
$$t(1,10,15,30;n)=\f 23N(1,10,15,30;8n+56)-2N(1,10,15,30;2n+14).$$
\endpro

\pro{Conjecture 2.18} Let $n\in\Bbb N$ with $n\e 2\mod 8$. Then
$$t(1,7,28,28;n)=\f 23N(1,7,28,28;8n+64)-2N(1,7,28,28;2n+16).$$
\endpro

\pro{Conjecture 2.19} Let $n\in\Bbb N$ with $n\e 8\mod 9$. Then
$$t(1,9,16,18;n)=\f 23N(1,9,16,18;8n+44)-2N(1,9,16,18;2n+11).$$
\endpro

\pro{Conjecture 2.20} Let $n\in\Bbb N$ with $n\e 1,7\mod 9$. Then
$$t(1,9,18,24;n)=\f 23N(1,9,18,24;8n+52)-2N(1,9,18,24;2n+13).$$
\endpro

\pro{Conjecture 2.21} Let $n\in\Bbb N$ with $n\e 1,4\mod 9$. Then
$$t(1,9,18,32;n)=\f 23N(1,9,18,32;8n+60)-2N(1,9,18,32;2n+15).$$
\endpro

\pro{Conjecture 2.22} Let $n\in\Bbb N$ with $n\e 5\mod 9$. Then
$$t(1,9,18,40;n)=\f 23N(1,9,18,40;8n+68)-2N(1,9,18,40;2n+17).$$
\endpro

\pro{Conjecture 2.23} Let $n\in\Bbb N$ with $n\e 2,5\mod 9$. Then
$$t(1,10,27,30;n)=\f 23N(1,10,27,30;8n+68)-2N(1,10,27,30;2n+17).$$
\endpro

\end{document}